\documentclass[12pt,leqno,twoside]{amsart}
\usepackage{amsmath}
\usepackage{amssymb}

\newtheorem{corollary}{Corollary}
\newtheorem{proposition}{Proposition}

\newenvironment{definition}
{\smallskip\noindent{\bf Definition\/}:}{\smallskip\par}

\newenvironment{remarks}
{\smallskip\noindent{\bf Remarks\/}.}{\smallskip\par}

\newcommand{\CC}{{\Bbb C}}

\newcommand{\ZZ}{{\Bbb Z}}

\newcommand{\DD}{{\Bbb D}}

\newcommand{\TT}{{\Bbb T}}
\newcommand{\mathcalO}{{\mathcal O}}

\newcommand{\mathcalJ}{{\mathcal J}}

\newcommand{\mathcalL}{{\mathcal L}}

\newcommand{\mathcalD}{{\mathcal D}}

\newcommand{\eps}{\varepsilon}

\newcommand{\codim}{{\rm codim\,}}

\title[Chern obstructions]{Chern obstructions for collections of 1-forms on singular varieties}
\author{\sc W.~Ebeling}

\address{Universit\"{a}t Hannover, Institut f\"{u}r Mathematik,
Postfach 6009, D-30060 Hannover, Germany.}

\email{ebeling@math.uni-hannover.de}

\author{\sc S.~M.~Gusein-Zade}

\address{Moscow State University, Faculty of Mechanics and Mathematics,
Moscow, 119992, Russia.}

\email{sabir@mccme.ru}

\thanks{Partially supported by the DFG-programme ''Global methods in
complex geometry'' (Eb 102/4--3), grants RFBR--04--01--00762,
NSh--1972.2003.1}

\dedicatory{Dedicated to Jean-Paul Brasselet on the occasion of his 60th birthday}



\date{}

\begin{document}

\maketitle

\begin{abstract}
We introduce a certain index of a collection of germs of 1-forms on a germ of a singular variety which
is a generalization of the local Euler obstruction corresponding to Chern numbers different from
the top one.
\end{abstract}

\section*{Introduction}

The aim of this paper is to bring together some ideas of \cite{BMPS} and \cite{EGLMS}.
A germ of a vector field or of a 1-form on the complex affine space $\CC^n$ at the
origin not vanishing in a punctured neighbourhood of it has a topological invariant~--- the Poincar\'e--Hopf
index. The sum of the Poincar\'e--Hopf indices of the singular points of a vector field on a compact complex
manifold is equal to the Euler characteristic of the manifold.
There are several generalizations of this notion to vector fields and/or to 1-forms
on complex analytic varieties with singularities (isolated or not) started by M.-H.~Schwartz:
\cite[\dots]{Schw, BS, GSV, BG, SS, EGMMJ, BMPS}. For the case of an isolated complete
intersection singularity  there is defined an index which is sometimes called the GSV index:
\cite{GSV, SS, EGMMJ}.  Another generalization which makes sense not only for isolated complete
intersection singularities and also not only for varieties with isolated singularities is the so
called local Euler obstruction:
\cite{BLS, BMPS} (its analogue for 1-forms was considered in \cite{EGGD}). One can say that
in some sense  all these invariants correspond to the Euler characteristic,
which, for a compact complex analytic manifold $M^n$, coincides with the top 
Chern number $\langle c_n(M), [M] \rangle$.

A generalization of the GSV-index corresponding to other Chern numbers (different from the
top one) was introduced and studied in \cite{EGLMS}. It is defined for a collection of germs of 1-forms
on an isolated complete intersection singularity. 
For a collection of 1-forms on a projective complex complete intersection with isolated singularities,
the sum of these indices of the singular points is
equal to plus-minus the corresponding Chern number of a smoothing of the variety. 

Here we define and
study an index of a collection of germs of 1-forms on a germ of a singular variety which is an analogue
of the local Euler obstruction corresponding to a Chern number different from the top one.

\section{Special points of 1-forms}
Let $(X^n,0)\subset(\CC^N,0)$ be the germ of a purely $n$-dimensional reduced
complex analytic variety at the origin (generally speaking with a non-isolated singularity).
Let ${\bf k}=\{k_i\}$, $i=1,\ldots, s$, be a fixed partition of $n$ (i.e., $k_i$ are positive
integers,
$\sum\limits_{i=1}^s k_i=n$). Let $\{\omega^{(i)}_j\}$ ($i=1,\ldots, s$,
$j=1,\ldots, n-k_i+1$) be a collection of germs of 1-forms on $(\CC^N, 0)$ (not necessarily
complex analytic; it suffices that the forms $\omega^{(i)}_j$ are complex linear functions
continuously depending on a point of $\CC^N$). Let $\eps>0$ be small enough so that there is
a representative $X$ of the germ $(X,0)$ and representatives $\omega^{(i)}_j$ of the germs of
1-forms inside the ball $B_\eps(0)\subset\CC^N$.

\begin{definition}
A point $P\in X$ is called a {\em special} point of the collection $\{\omega^{(i)}_j\}$  of
1-forms on the variety $X$ if there exists a sequence $\{P_m\}$ of points from the
non-singular part $X_{\rm reg}$ of the variety $X$ converging to $P$ such that the sequence $T_{P_m}X_{\rm
reg}$ of the tangent spaces at the points $P_m$ has a limit $L$ as $m \to \infty$ 
(in the Grassmann manifold of $n$-dimensional vector subspaces of $\CC^N$) and the
restrictions of the 1-forms $\omega^{(i)}_1$, \dots, $\omega^{(i)}_{n-k_i+1}$ to the subspace
$L\subset T_P\CC^N$ are linearly dependent for each $i=1, \ldots, s$.
\end{definition}

\begin{definition}
The collection $\{\omega^{(i)}_j\}$ of 1-forms has an {\em isolated special point} on the germ $(X,0)$
if it has no special points on $X$ in a punctured neighbourhood of the origin.
\end{definition}

\begin{remarks} {\bf 1.} 
If the 1-forms $\omega^{(i)}_j$ are complex analytic, the property to have an isolated
special point is a condition on the classes of these 1-forms in the module
$$
\Omega^1_{X,0}=
\Omega^1_{\CC^N,0}/\{f\cdot\Omega^1_{\CC^N,0} + df\cdot\mathcalO_{\CC^N,0}\vert f\in\mathcalJ_X\}
$$
of germs of 1-forms on the variety $X$ ($\mathcalJ_X$ is the ideal of germs of
holomorphic functions vanishing on $X$).

\noindent
{\bf 2.} For the case $s=1$ (and therefore $k_1=n$), i.e. for one 1-form $\omega$, there
exists a notion of a {\em singular} point of the 1-form $\omega$ on $X$ (see, e.g.,
\cite{EGGD}). It is defined in terms of a Whitney stratification of the variety $X$. A point $x
\in X$ is a {\em singular} point of the 1-form $\omega$ on the variety $X$ if the restriction of the
1-form $\omega$ to the stratum of $X$ containing $x$ is equal to zero at the point $x$. (One
should consider points of all zero-dimensional strata as singular ones.)  One can easily see
that a special point of the 1-form $\omega$ on the variety $X$ is singular, but not vice versa.
(E.g.\ the origin is a singular point of the 1-form $dx$ on the cone $\{x^2+y^2+z^2=0\}$, but
not a special one.) On a smooth variety these two notions coincide.
\end{remarks}

The notion of a {\em non-degenerate} special (singular) point of a collection of germs of 1-forms
on a  smooth variety was introduced in \cite{EGLMS}. The index of a non-degenerate point of a
collection of germs of holomorphic 1-forms is equal to 1.

Let 
$$
\mathcalL^{\bf k}= \prod\limits_{i=1}^s \prod\limits_{j=1}^{n-k_i+1} \CC^{N\ast}_{ij}
$$
be the space of collections of linear functions on $\CC^N$ (i.e.\ of 1-forms with constant
coefficients).

\begin{proposition}\label{prop1}
There exists an open and dense subset $U \subset \mathcalL^{\bf k}$ such that each collection
$\{\ell^{(i)}_j\} \in U$ has only isolated special points on $X$ and, moreover, all these points
belong to the smooth part $X_{\rm reg}$ of the variety $X$ and are non-degenerate.
\end{proposition}

\begin{proof}
Let $Y\subset X\times \mathcalL^{\bf k}$ be the closure of the set of pairs $(x, \{\ell^{(i)}_j\})$
where $x\in X_{\rm reg}$ and the restrictions of the linear functions $\ell^{(i)}_1$, \dots,
$\ell^{(i)}_{n-k_i+1}$ to the tangent space $T_xX_{\rm reg}$ are linearly dependent for each
$i=1, \ldots, s$. Let $\pi:Y\to \mathcalL^{\bf k}$ be the projection to the second factor.
One has $\codim Y= \sum\limits_{i=1}^s k_i=n$ and therefore $\dim Y = \dim \mathcalL^{\bf k}$.
Moreover, $Y\setminus (X_{\rm reg}\times \mathcalL^{\bf k})$ is a proper subvariety of $Y$
and therefore its dimension is strictly smaller than $\dim \mathcalL^{\bf k}$. A generic point of
the space $\mathcalL^{\bf k}$ is a regular value of the map $\pi$ which means that it has only
finitely many preimages, all of them belong to $X_{\rm reg}\times \mathcalL^{\bf k}$ and the map $\pi$
is non-degenerate at them. This implies the statement.
\end{proof}

\begin{corollary} \label{cor1}
Let $\{\omega^{(i)}_j\}$ be a collection of 1-forms on $X$ with an isolated special
point at the origin. Then there exists a deformation $\{\widetilde \omega^{(i)}_j\}$
of the collection $\{\omega^{(i)}_j\}$ whose special points lie in $X_{\rm reg}$ and are
non-degenerate. Moreover, as such a deformation one can use $\{\omega^{(i)}_j+
\lambda \ell^{(i)}_j\}$ with a generic collection $\{\ell^{(i)}_j\}\in \mathcalL^{\bf k}$, $\lambda \neq 0$
small enough.
\end{corollary}

\begin{corollary} \label{cor2}
The set of collections of holomorphic 1-forms with a non-isolated special point at the origin
has infinite codimension in the space of all holomorphic collections.
\end{corollary}


\section{Local Chern obstructions}
Let $\{\omega^{(i)}_j\}$ be a collection of germs of 1-forms on $(X, 0)$ with an isolated special
point at the origin. Let $\nu : \widehat{X} \to X$ be the Nash
transformation of the variety $X\subset B_\eps(0)$ defined as follows. Let $G(n,N)$ be the
Grassmann manifold
of $n$-dimensional vector subspaces of $\CC^N$. There is a natural map
$\sigma : X_{\rm reg} \to B_\eps(0) \times G(n,N)$ which
sends a point $x\in X_{\rm reg}$ to $(x, T_x X_{\rm reg})$.
The Nash transform $\widehat{X}$ of the variety $X$ is the closure of the image ${\rm Im}\,
\sigma$ of the map $\sigma$ in $B_\eps(0) \times G(n,N)$, $\nu$ is the natural projection.
The Nash bundle $\widehat{T}$ over $\widehat{X}$ is a
vector bundle of rank $n$
which is the pullback of the tautological bundle on the Grassmann manifold
$G(n,N)$. There is a
natural lifting of the Nash transformation to a bundle map from the Nash bundle
$\widehat{T}$ to the restriction of the tangent bundle $T\CC^N$ of $\CC^N$
to $X$. This is an isomorphism of $\widehat{T}$ and $TX_{\rm reg} \subset T\CC^N$
over the non-singular part $X_{\rm reg}$ of $X$.

The collection of 1-forms $\{\omega^{(i)}_j\}$ gives rise to a
section $\widehat{\omega}$ of the bundle
$$
\widehat\TT=\bigoplus_{i=1}^s\bigoplus_{j=1}^{n-k_i+1}\widehat T^*_{i,j}
$$
where $\widehat{T}^\ast_{i,j}$ are copies of the dual Nash bundle $\widehat{T}^\ast$ over
the Nash transform $\widehat{X}$ numbered by indices $i$ and $j$. Let
$\widehat\DD\subset\widehat\TT$ be the set of pairs $(x,\{\alpha^{(i)}_j\})$ where $x\in\widehat
X$ and the collection $\{\alpha^{(i)}_j\}$ of elements of $\widehat T_x^*$ (i.e.\ of linear
functions on
$\widehat T_x$) is such that $\alpha^{(i)}_1$, \dots, $\alpha^{(i)}_{n-k_i+1}$ are linearly
dependent for each $i=1, \dots, s$. The image of the section $\widehat\omega$ does not intersect
$\widehat\DD$ outside of the preimage $\nu^{-1}(0)\subset\widehat X$ of the origin. The
map $\widehat\TT\setminus \widehat\DD\to \widehat X$ is a fibre bundle. The fibre
$W_x=\widehat\TT_x \setminus \widehat \DD_x$ of it is $(2n-2)$-connected, its homology group
$H_{2n-1}(W_x;\ZZ)$ is isomorphic to $\ZZ$ and has a natural generator:
see, e.g., \cite{EGLMS}. The latter fact implies that the fibre bundle $\widehat\TT\setminus
\widehat\DD\to \widehat X$ is homotopically simple in dimension $2n-1$, i.e.\ the fundamental
group $\pi_1(\widehat X)$ of the base acts trivially on the homotopy group
$\pi_{2n-1}(W_x)$ of the fibre, the last one being isomorphic to the homology
group $H_{2n-1}(W_x;\ZZ)$: see, e.g., \cite{St}. 

\begin{definition}
The {\em local Chern obstruction} ${\rm Ch}_{X,0}\,\{\omega^{(i)}_j\}$ of the collections
of germs of 1-forms $\{\omega^{(i)}_j\}$ on $(X,0)$ at the origin is the (primary, and in fact
the only) obstruction to extend the section $\widehat{\omega}$ of the fibre bundle
$\widehat\TT\setminus
\widehat\DD\to
\widehat X$ from the preimage of a neighbourhood of the sphere $S_\eps= \partial B_\eps$ to 
$\widehat X$, more precisely its value (as an element of
$H^{2n}(\nu^{-1}(X\cap B_\eps), \nu^{-1}(X \cap S_\eps);\ZZ)$\,) on the
fundamental class of the pair $(\nu^{-1}(X\cap B_\eps), \nu^{-1}(X \cap S_\eps))$.
\end{definition}

The definition of the local Chern obstruction ${\rm Ch}_{X,0}\,\{\omega^{(i)}_j\}$ 
can be reformulated in the following way. Let $\mathcalD^{\bf k}_X\subset \CC^N \times \mathcalL^{\bf k}$
be the closure of the set of pairs $(x, \{\ell^{(i)}_j\})$ such that $x \in X_{\rm reg}$ and the
restrictions of the linear functions $\ell_1^{(i)}$, \dots , $\ell_{n-k_i+1}^{(i)}$ to
$T_xX_{\rm reg} \subset \CC^N$ are linearly dependent for each $i=1, \ldots, s$. (For $s=1$,
${\bf k}=\{n\}$, $\mathcalD^{\bf k}_X$ is the (non-projectivized) conormal space of $X$ \cite{Te}.)
The collection $\{\omega^{(i)}_j\}$ of germs of 1-forms on $(\CC^N,0)$ defines a section
$\check{\omega}$ of the (trivial) fibre bundle $\CC^N \times \mathcalL^{\bf k} \to \CC^N$. Then
$$
{\rm Ch}_{X,0}\,\{\omega^{(i)}_j\} = ( \check{\omega}(\CC^N) \circ \mathcalD^{\bf k}_X)_0
$$
where $(\cdot \circ \cdot )_0$ is the intersection number at the origin in $\CC^N \times
\mathcalL^{\bf k}$. This description can be considered as a generalization of an expression of the
local Euler obstruction as a microlocal intersection number defined in \cite{KS}, see also
\cite[Sections 5.0.3 and 5.2.1]{Schu}.

\begin{sloppypar}

\begin{remarks} {\bf 1.} On a smooth manifold $X$ the local Chern obstruction
${\rm Ch}_{X,0}\,\{\omega^{(i)}_j\}$ coincides with the index ${\rm ind}_{X,0} \,
\{\omega^{(i)}_j\}$ of the collection $\{\omega^{(i)}_j\}$ defined in \cite{EGLMS}. 

\noindent
{\bf 2.} The local Euler obstruction is defined for vector fields as well as for 1-forms. One
can see that vector fields are not well adapted to a definition of the local Chern obstruction.
A more or less direct version of the definition above for vector fields demands to consider
vector fields on a singular variety $X \subset \CC^N$ to be sections $v=v(x)$ of 
$T\CC^N\mbox{\raisebox{-0.2ex}{$\vert$}}{}_{X}$
such that $v(x) \in T_xX \subset T_x\CC^N$ ($\dim T_xX$ is not constant). (Traditionally vector
fields tangent to smooth strata of the variety $X$ are considered.) There exist only continuous
(non-trivial, i.e.\ with $s>1$) collections of such vector fields ''on $X$'' with isolated
special points, but not holomorphic ones.

\noindent
{\bf 3.} The definition of the local Chern obstruction ${\rm Ch}_{X,0}\,\{\omega^{(i)}_j\}$ may
also be formulated in terms of a collection $\{\omega^{(i)}\}$ of germs of 1-forms with values in
vector spaces $L_i$ of dimensions $n-k_i+1$. Therefore (via differentials) it is also defined
for a collection $\{f^{(i)}\}$ of germs of maps $f^{(i)} : (\CC^N,0) \to (\CC^{n-k_i+1},0)$
(just as the Euler obstruction is defined for a germ of a function).
\end{remarks}

\end{sloppypar}

Being a (primary) obstruction, the local Chern obstruction satisfies the law of conservation of
number, i.e.\ if a collection of 1-forms $\{\widetilde\omega^{(i)}_j\}$ is a deformation of the
collection $\{\omega^{(i)}_j\}$ and has isolated special points on $X$, then 
$$
{\rm Ch}_{X,0} \, \{\omega^{(i)}_j\} = \sum{\rm Ch}_{X,Q}\, \{\widetilde\omega^{(i)}_j\}
$$
where the sum on the right hand side is over all special points $Q$
of the collection $\{\widetilde\omega^{(i)}_j\}$ on $X$ in a neighbourhood of the origin.
With Corollary~\ref{cor1} this implies the following statements.

\begin{proposition}\label{prop2}
The local Chern obstruction ${\rm Ch}_{X,0} \, \{\omega^{(i)}_j\}$ of a collection
$\{\omega^{(i)}_j\}$ of germs of holomorphic 1-forms is equal to the number of special points on
$X$ of a generic (holomorphic) deformation of the collection.
\end{proposition}

This statement is an analogue of Proposition~2.3 in \cite{STV}.

\begin{proposition}\label{prop3}
If a collection $\{\omega^{(i)}_j\}$ of 1-forms on a compact (say, projective) variety $X$ has
only isolated special points, then the sum of the local Chern obstructions of the
collection $\{\omega^{(i)}_j\}$ at these points does not depend on the collection and
therefore is an invariant of the variety.
\end{proposition}

It is reasonable to consider this sum as ($(-1)^n$ times) the corresponding Chern number of the
singular variety $X$.

Let $(X,0)$ be an isolated complete intersection singularity. As it was mentioned above, a
collection of germs of 1-forms $\{\omega^{(i)}_j\}$ on $(X,0)$ with an isolated special point at the
origin has an index ${\rm ind}_{X,0}\, \{\omega^{(i)}_j\}$ which is an analogue of the GSV--index
of a vector field: \cite{EGMMJ}. The fact that both the Chern obstruction and the index satisfy the law of
conservation of number and they coincide on a smooth manifold yields the following statement.

\begin{proposition}\label{prop4}
For a collection $\{\omega^{(i)}_j\}$ of germs of 1-forms on an isolated complete intersection singularity
$(X,0)$ the difference
$$
{\rm ind}_{X,0}\, \{\omega^{(i)}_j\} - {\rm Ch}_{X,0} \, \{\omega^{(i)}_j\}
$$
does not depend on the collection and therefore is an invariant of the germ of the variety.
\end{proposition}

Since, by Proposition~\ref{prop1}, ${\rm Ch}_{X,0} \, \{\ell^{(i)}_j\}=0$ for a generic collection
$\{\ell^{(i)}_j\}$ of linear functions on $\CC^N$, one has the following statement.

\begin{corollary}
One has
$$
{\rm Ch}_{X,0} \, \{\omega^{(i)}_j\} = {\rm ind}_{X,0}\, \{\omega^{(i)}_j\} - {\rm ind}_{X,0}\,
\{\ell^{(i)}_j\}
$$
for a generic collection $\{\ell^{(i)}_j\}$ of linear functions on $\CC^N$.
\end{corollary}


\begin{thebibliography}{12}


\bibitem{BG} Ch.~Bonatti, X.~G\'omez-Mont: The index of holomorphic
vector fields on singular varieties I. Ast\'erisque {\bf 222}, 9--35 (1994).

\bibitem{BLS} J.-P.~Brasselet, L\^e D\~ung Tr\'ang, J.~Seade: Euler obstruction
and indices of
vector fields. Topology {\bf 39}, 1193--1208 (2000).

\bibitem{BMPS} J.-P.~Brasselet, D.~Massey, A.~J.~Parameswaran,
J.~Seade: Euler obstruction and defects of functions on singular
varieties. J. London Math. Soc. (2) {\bf 70}, 59--76 (2004).

\bibitem{BS} J.-P.~Brasselet, M.-H.~Schwartz: Sur les classes de Chern
d'un ensemble analytique complexe. In:
Caract\'eristique d'Euler-Poincar\'e,  Ast\'erisque {\bf 82--83}, 93--147 (1981).

\bibitem{EGMMJ} W.~Ebeling, S.~M.~Gusein-Zade: Indices of 1-forms
on an isolated complete intersection singularity. Moscow Math. J.
{\bf 3}, 439--455 (2003).

\bibitem{EGLMS} W.~Ebeling, S.~M.~Gusein-Zade: Indices of vector fields or 1-forms and
characteristic numbers. math.AG/0303330. Bull. London Math. Soc. (to appear).

\bibitem{EGGD} W.~Ebeling, S.~M.~Gusein-Zade: Radial index and Euler obstruction of a
1-form on a singular variety. math.AG/0402388. Geom. Dedicata (to appear).

\bibitem{GSV} X.~G\'omez-Mont, J.~Seade,
A.~Verjovsky: The index of a holomorphic flow with an isolated singularity. Math. Ann.
{\bf 291}, 737--751 (1991).

\bibitem{KS} M.~Kashiwara, P.~Schapira: Sheaves on Manifolds.
Springer-Verlag, 1990.

\bibitem{Schu} J.~Sch\"urmann: Topology of Singular Spaces and
Constructible Sheaves. Birkh\"auser, 2003.

\bibitem{Schw} M.-H.~Schwartz: Classes caract\'eristiques
d\'efinies par une stratification d'une vari\'et\'e analytique complexe.
C. R. Acad. Sci. Paris S\'er. I Math. {\bf 260}, 3262--3264, 3535--3537 (1965).

\bibitem{SS} J.~A.~Seade, T.~Suwa: A residue formula for the
index of a holomorphic flow. Math. Ann. {\bf 304}, 621--634 (1996).

\bibitem{STV} J.~Seade, M.~Tib\v{a}r, A.~Verjovsky: Milnor numbers and Euler obstruction.
math.CV/0311109.

\bibitem{St} N.~Steenrod: The Topology of Fibre Bundles. Princeton
Math. Series, Vol.~{\bf 14}, Princeton University Press, Princeton, N. J., 1951.

\bibitem{Te} B.~Teissier:
Vari\'et\'es polaires. II. Multiplicit\'es polaires, sections planes, et conditions de Whitney. 
In: Algebraic geometry (La R\'abida, 1981), 
Lecture Notes in Math., Vol.~{\bf 961},
Springer, Berlin, 1982, pp.~314--491. 

\end{thebibliography}
\end{document}